\documentclass[11pt]{article}
\usepackage{amssymb}
\usepackage{amsmath}
\usepackage{amsthm}
\usepackage{enumerate}
\usepackage[usenames]{color}
\usepackage[dvipsnames]{xcolor}
\usepackage{graphicx}
\usepackage{subfigure}
\usepackage{enumerate}
\usepackage[usenames]{color}
\definecolor{brown}{cmyk}{0, 0.72, 1, 0.45}
\definecolor{grey}{gray}{0.5}

\usepackage[normalem]{ulem}

\usepackage{tikz}
\usetikzlibrary{arrows}
\usetikzlibrary{decorations}
\usetikzlibrary{shapes.misc}
\newcommand{\pic}[1]{
\begin{tikzpicture}
#1
\end{tikzpicture}}
\newcommand{\pics}[2]{
\begin{tikzpicture}[scale=#1]
#2
\end{tikzpicture}}

\newcommand{\lem}[3]{\begin{lemma}\label{#1}#2\end{lemma}\proofstart #3\proofend}

\allowdisplaybreaks
\title{}

\def\a{\alpha} \def\b{\beta}  
    
\def\G{\Gamma}  
     
   \def\p{\pi}
  \def\s{\sigma}

\def\cP{{\cal P}}

\newtheoremstyle{plain}%
{8pt plus2pt minus4pt}%
{8pt plus2pt minus4pt}%
{\itshape}%
{}%
{\bfseries\scshape}%
{}%
{6pt}
{}%

\newtheoremstyle{remark}%
{8pt plus2pt minus4pt}%
{8pt plus2pt minus4pt}%
{\upshape}
{}%
{\bfseries\scshape}%
{}%
{6pt}
{}%

\theoremstyle{plain}
\newtheorem{theorem}{Theorem}
\newtheorem{lemma}[theorem]{Lemma}

\theoremstyle{remark}

\newcommand{\proofstart}{{\bf Proof\hspace{2em}}}
\newcommand{\proofend}{\hspace*{\fill}\mbox{$\Box$}}

\newcommand{\beq}[2]{\begin{equation}\label{#1}#2\end{equation}}
\newcommand{\eeq}{\end{equation}}

\newcommand{\set}[1]{\left\{#1\right\}}

\def\Pr{\mbox{{\bf Pr}}}

\newcommand{\ignore}[1]{}

\def\a{\alpha} \def\b{\beta}  
    
\def\G{\Gamma}  
     
   \def\p{\pi}
  \def\s{\sigma}

\def\cP{{\cal P}}

\def\Pr{\mbox{{\bf Pr}}}

\begin{document}

\title{Adding random edges to create the square of a Hamilton cycle}
\author{Patrick Bennett\thanks{Department of Mathematics, Department of Mathematics, Western Michigan University, Kalamazoo, MI. Supported in part by Simons Foundation Grant \#426894.}\and Andrzej Dudek\thanks{Department of Mathematics, Department of Mathematics, Western Michigan University, Kalamazoo, MI. Supported in part by Simons Foundation Grant \#522400.}\and Alan Frieze\thanks{Department of Mathematical Sciences, Carnegie Mellon University, Pittsburgh PA. Research supported in part by NSF grant CCF1013110 and Simons Foundation Grant \#333329.}}

\date{\today}

\maketitle \makeatother

\begin{abstract}
We consider how many random edges need to be added to a graph of order $n$ with minimum degree $\a n$ in order that it contains the square of a Hamilton cycle w.h.p..
\end{abstract}

\section{Introduction}

By the \emph{$k$th power of a Hamilton cycle}, we mean a permutation (bijection) $\pi:[n] \rightarrow [n]$ such that $\set{\pi(i),\pi(j)}\in E(G)$ whenever $i<j\leq i+k$. (Here $i+k$ is to be taken as $i+k-n$ if $i+k\geq n+1$.) Hamilton cycles have long been studied in the context of random graphs (see, e.g.,  \cite{AKS, Boll, KS, Posa}). Powers of Hamilton cycles are less well-studied and much less is known about them.  

K\"uhn and Osthus \cite{KO} observed that for $k\geq 3$, $p=\frac{1}{n^{1/k}}$ is the coarse threshold for the existence of the $k$th power of a Hamilton cycle in $G_{n,p}$. This comes directly from a result of Riordan~\cite{Rio}. For $k=2$ they gave a bound of $p\geq n^{-1/2+\varepsilon}$ (for any $\varepsilon>0$) being sufficient for the existence of the square of a Hamilton cycle w.h.p.. This result was improved by Nenadov and \v{S}kori\'c~\cite{NS} to $p\geq \frac{C\log^4n}{\sqrt{n}}$ ($C$ is a positive constant) being sufficient for the existence of the square of a Hamilton cycle. 

In this paper we consider a problem related to the Pos\'a-Seymour conjecture, which states that every graph $G$ on $n$ vertices with minimum degree at least $kn/(k+1)$ contains the $k$th power of a Hamilton cycle. This conjecture was proved for large enough $n$ by Koml\'os, Sark\"ozy and Szemer\'edi \cite{KSS}. Bohman, Frieze and Martin \cite{BFM} considered the question of how many random edges need to be added to a graph with minimum degree $\a n$ with $0<\a<1/2$ in order that it is Hamiltonian w.h.p.. They showed that $(30\log\a^{-1}+13)n$ random edges are sufficient.

The following theorem extends their result to the square of a Hamilton cycle. For graphs $G=(V,E)$ and $X=(V,F)$ we define a graph $G+X$ on vertex set $V$ with edge set $E\cup F$.
\begin{theorem}\label{th1}
For every constant $\a>1/2$ there exists a sufficiently large $K=K(\a)$ such that the following holds. Let $G$ be a graph of order $n$ which has minimum degree at least $\a n$. Let $X$ denote a set of randomly chosen edges. Then w.h.p. $\G=G+X$ contains the square of a Hamilton cycle, provided that
$$|X|\geq Kn^{4/3}\log^{1/3}n.$$
\end{theorem}
Clearly $n^{4/3}$ is less than the $n^{3/2}$ needed if all edges are random.
\section{Proof of Theorem \ref{th1}}

\subsection{Preliminaries}

It will be convenient for the computations to assume that the edges $X$ will be given as $X=X_1\cup X_2\cup X_3$ where each of the sets in this partition are independent random subsets of $E(K_n)$ where each edge is independently included with probability $p= \frac{K\log^{1/3}n}{n^{2/3}}$ for some large positive constant~$K$.

Assume first that $n=2m$ is even. It follows from Erd\H{o}s and R\'enyi \cite{ER1} that w.h.p. the edges $X_1$ contain a perfect matching $M$. By symmetry, it will be a random matching of $K_n$ that is independent of $G$. It can therefore be derived from a random permutation $\p=(z_1,z_2,\ldots,z_n)$ via $M=\set{e_1,e_2,\ldots,e_m}$ where $e_i=\set{z_{2i-1},z_{2i}},i=1,2,\ldots,m$.

Now for any graph $H$ on vertex set $V(G)$, define a graph $\Pi(H)$ with vertex set $M$ and an edge $\set{e,f}, e,f\in M$ whenever the subgraph $H_{e,f}$ of $H$ induced by the four vertices in $e\cup f$ is $K_4$. Let $\G_1=\Pi(G+X_1)$.
We argue that 
\beq{mindeg}{
\text{w.h.p., }\G_1 \text{ has minimum degree at least $\b_1n$ where $\b_1=(2\a-1)^3/2$.}
}
To see this consider a fixed edge $e=\set{x,y}\in \G_1$. Let $N(a)$ denote the set of neighbors of vertex $a$ in $G$ and note that $|N(x)\cap N(y)|\geq (2\a-1)n$. The probability that another edge $\set{u,v}\in \G_1$ satisfies $u,v\in N(x)\cap N(y)$ is at least $(1-o(1))(2\a-1)^2$. Thus the degree $d_e$ of edge $e$ has expectation at least $(1-o(1))(2\a-1)^3n$ in $\G_1$. Swapping a pair in permutation $\p$ can only change $d_e$ by at most one. Applying a version of the Azuma-Hoeffding inequality (see for example McDiarmid \cite{McD} or Frieze and Pittel \cite{FP}) we see that $\Pr(d_e\leq (2\a-1)^3n/2)\leq e^{-\Omega((2\a-1)^6n)}$. This verifies \eqref{mindeg}, after inflating the probability bound by $m$. Note that only the edges of $\G_1$ are used here. 

Now let $\G_2=\Pi(G+X_1 + X_2)$.
\begin{lemma}\label{Gamcon}
$\G_2$ is connected w.h.p.. 
\end{lemma}
\proofstart
It follows from \eqref{mindeg} that $\G_1$ has at most $1/\b_1$ components and these are all of size at least $\b_1n$. Thus,
\[
\Pr(\G_2\text{ is not connected})\leq\frac{1}{\b_1}\max_{\b_1n\leq s\leq m/2} (1-p^3)^{s(\a n-2s)/2}=o(1).
\]
Indeed, fix a component $S$ of $\G_1$. It will have size $s\in[\b_1n/2,m/2]$. For each $e=\{u,v\}\in S$ there are at least $(\a n-2s)$ vertices $T$ outside $\bigcup_{e\in S}e$. For each vertex $x\in T$ we have a matching edge $f=\{x,y\}\in \bar{S}$ such that $e$ and $f$ are joined by an edge (say $\{u,x\}$) from $G$. The term $p^3$ accounts for the probability that $X_2$ will provide another three edges ($\{u,y\}$, $\{v,x\}$ and $\{v,y\}$) to create a $K_4$. We divide by two in $s(\a n-2s)/2$ to account for there being two choices for $x\in f$. 
\proofend
\subsection{2-paths}\label{2path}
\begin{figure}[h]
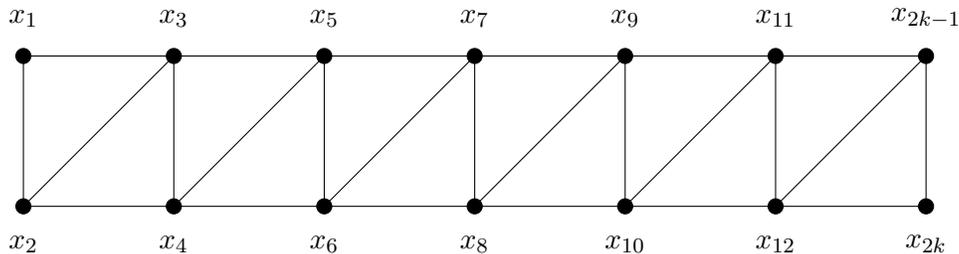

\begin{center}
\pic{
\draw (1,0) -- (13,0);
\draw (1,-2) -- (13,-2);
\draw (1,0) -- (1,-2);
\draw (3,0) -- (3,-2);
\draw (5,0) -- (5,-2);
\draw (7,0) -- (7,-2);
\draw (9,0) -- (9,-2);
\draw (11,0) -- (11,-2);
\draw (13,0) -- (13,-2);
\draw (1,-2) -- (3,0);
\draw (3,-2) -- (5,0);
\draw (5,-2) -- (7,0);
\draw (7,-2) -- (9,0);
\draw (9,-2) -- (11,0);
\draw (11,-2) -- (13,0);
\node at (1,.5) {$x_1$};
\node at (3,.5) {$x_3$};
\node at (5,.5) {$x_5$};
\node at (7,.5) {$x_7$};
\node at (9,.5) {$x_9$};
\node at (11,.5) {$x_{11}$};
\node at (13,.5) {$x_{2k-1}$};
\node at (1,-2.5) {$x_2$};
\node at (3,-2.5) {$x_4$};
\node at (5,-2.5) {$x_6$};
\node at (7,-2.5) {$x_8$};
\node at (9,-2.5) {$x_{10}$};
\node at (11,-2.5) {$x_{12}$};
\node at (13,-2.5) {$x_{2k}$};

\draw [fill=black] (1,0) circle [radius=0.1];
\draw [fill=black] (3,0) circle [radius=0.1];
\draw [fill=black] (5,0) circle [radius=0.1];
\draw [fill=black] (7,0) circle [radius=0.1];
\draw [fill=black] (9,0) circle [radius=0.1];
\draw [fill=black] (11,0) circle [radius=0.1];
\draw [fill=black] (13,0) circle [radius=0.1];

\draw [fill=black] (1,-2) circle [radius=0.1];
\draw [fill=black] (3,-2) circle [radius=0.1];
\draw [fill=black] (5,-2) circle [radius=0.1];
\draw [fill=black] (7,-2) circle [radius=0.1];
\draw [fill=black] (9,-2) circle [radius=0.1];
\draw [fill=black] (11,-2) circle [radius=0.1];
\draw [fill=black] (13,-2) circle [radius=0.1];

}
\end{center}
\caption{A 2-path for $k=7$.}
\label{fig:2path}
\end{figure}

A {\em 2-path} is a sequence of vertices $(x_1,x_2,\ldots,x_{2k})$ such that (i) $(x_1,x_2,\ldots,x_{2k})$ is a path in $G+X$, (ii) $\set{x_{2i-1},x_{2i}}\in M,\,i=1,2,\ldots,k$, and (iii) $\set{x_i,x_{i+2}}$ are edges of $\G=G+X$ for $i=1,2,\ldots,2k-2$ (see Figure~\ref{fig:2path}).

In a 2-path we refer to the edges $\set{x_{2i-1},x_{2i}},i=1,2,\ldots,k$ as the {\em pillars}.

\begin{figure}
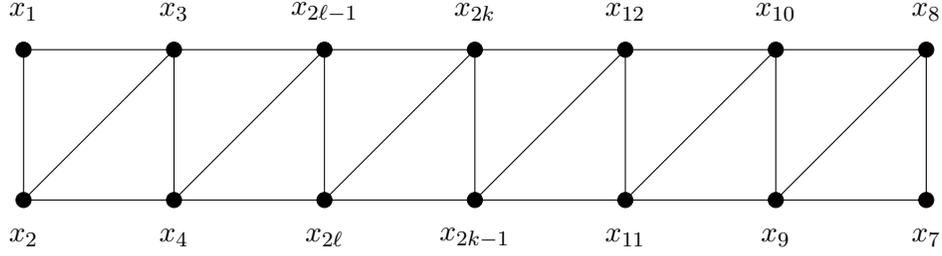

\begin{center}
\pic{
\draw (1,0) -- (13,0);
\draw (1,-2) -- (13,-2);
\draw (1,0) -- (1,-2);
\draw (3,0) -- (3,-2);
\draw (5,0) -- (5,-2);
\draw (7,0) -- (7,-2);
\draw (9,0) -- (9,-2);
\draw (11,0) -- (11,-2);
\draw (13,0) -- (13,-2);
\draw (1,-2) -- (3,0);
\draw (3,-2) -- (5,0);
\draw (5,-2) -- (7,0);
\draw (7,-2) -- (9,0);
\draw (9,-2) -- (11,0);
\draw (11,-2) -- (13,0);
\node at (1,.5) {$x_1$};
\node at (3,.5) {$x_3$};
\node at (5,.5) {$x_{2\ell-1}$};
\node at (7,.5) {$x_{2k}$};
\node at (9,.5) {$x_{12}$};
\node at (11,.5) {$x_{10}$};
\node at (13,.5) {$x_{8}$};
\node at (1,-2.5) {$x_2$};
\node at (3,-2.5) {$x_4$};
\node at (5,-2.5) {$x_{2\ell}$};
\node at (7,-2.5) {$x_{2k-1}$};
\node at (9,-2.5) {$x_{11}$};
\node at (11,-2.5) {$x_{9}$};
\node at (13,-2.5) {$x_{7}$};

\draw [fill=black] (1,0) circle [radius=0.1];
\draw [fill=black] (3,0) circle [radius=0.1];
\draw [fill=black] (5,0) circle [radius=0.1];
\draw [fill=black] (7,0) circle [radius=0.1];
\draw [fill=black] (9,0) circle [radius=0.1];
\draw [fill=black] (11,0) circle [radius=0.1];
\draw [fill=black] (13,0) circle [radius=0.1];

\draw [fill=black] (1,-2) circle [radius=0.1];
\draw [fill=black] (3,-2) circle [radius=0.1];
\draw [fill=black] (5,-2) circle [radius=0.1];
\draw [fill=black] (7,-2) circle [radius=0.1];
\draw [fill=black] (9,-2) circle [radius=0.1];
\draw [fill=black] (11,-2) circle [radius=0.1];
\draw [fill=black] (13,-2) circle [radius=0.1];
}
\end{center}
\caption{A rotated path for $k=7$ and $\ell=3$.}
\label{fig:2pathrotation}
\end{figure}

We now define a {\em rotation} with $\set{x_1,x_2}$ as the {\em fixed end} and $\set{x_{2k-1},x_{2k}}$ as the {\em rotated end}. Suppose that for some $\ell\leq k-2$ we have that $\set{x_{2\ell-1},x_{2k}},\set{x_{2\ell},x_{2k-1}}, \set{x_{2\ell},x_{2k}}$ are all edges of $\G$. Then we obtain a new 2-path 
$(x_1,x_2,\ldots,x_{2\ell},x_{2k},x_{2k-1},\ldots,x_{2\ell+1})$ (see Figure~\ref{fig:2pathrotation}). 


\subsection{Algorithm {\sc ERA}}
\begin{figure}[h]
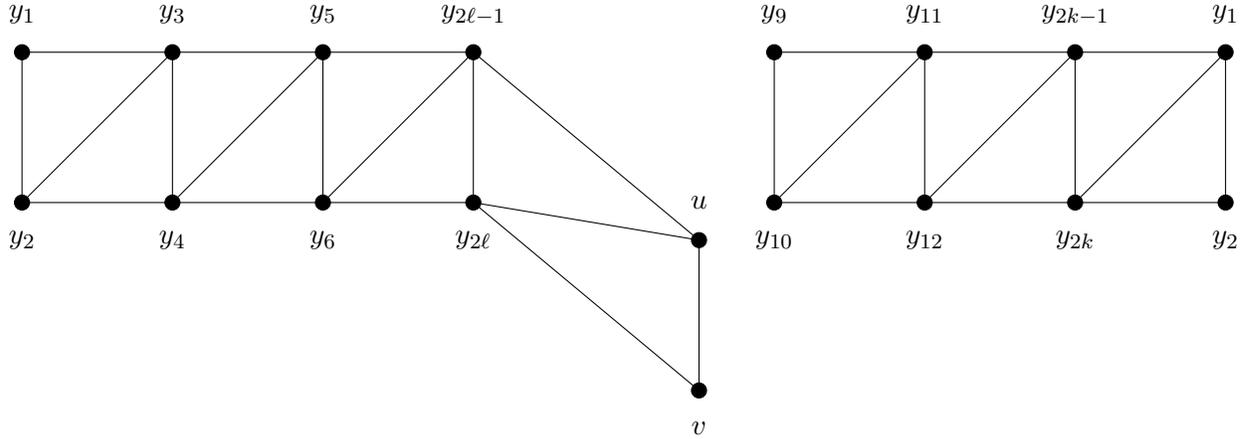

\begin{center}
\pic{
\draw (1,0) -- (7,0);
\draw (11,0) -- (17,0);
\draw (1,-2) -- (7,-2);
\draw (11,-2) -- (17,-2);
\draw (1,0) -- (1,-2);
\draw (3,0) -- (3,-2);
\draw (5,0) -- (5,-2);
\draw (7,0) -- (7,-2);
\draw (17,0) -- (17,-2);
\draw (11,0) -- (11,-2);
\draw (13,0) -- (13,-2);
\draw (15,0) -- (15,-2);
\draw (1,-2) -- (3,0);
\draw (3,-2) -- (5,0);
\draw (5,-2) -- (7,0);
\draw (11,-2) -- (13,0);
\draw (13,-2) -- (15,0);
\draw (15,-2) -- (17,0);
\node at (1,.5) {$y_1$};
\node at (3,.5) {$y_3$};
\node at (5,.5) {$y_5$};
\node at (7,.5) {$y_{2\ell-1}$};
\node at (11,.5) {$y_9$};
\node at (13,.5) {$y_{11}$};
\node at (15,.5) {$y_{2k-1}$};
\node at (1,-2.5) {$y_2$};
\node at (3,-2.5) {$y_4$};
\node at (5,-2.5) {$y_6$};
\node at (7,-2.5) {$y_{2\ell}$};
\node at (11,-2.5) {$y_{10}$};
\node at (13,-2.5) {$y_{12}$};
\node at (15,-2.5) {$y_{2k}$};
\node at (17,.5) {$y_1$};
\node at (17,-2.5) {$y_2$};
\draw (7,-2) -- (10,-4.5);
\draw (7,0) -- (10,-2.5);
\draw (10,-4.5) -- (10,-2.5);
\draw (7,-2) -- (10,-2.5);
\node at (10,-5) {$v$};
\node at (10,-2) {$u$};

\draw [fill=black] (1,0) circle [radius=0.1];
\draw [fill=black] (3,0) circle [radius=0.1];
\draw [fill=black] (5,0) circle [radius=0.1];
\draw [fill=black] (7,0) circle [radius=0.1];
\draw [fill=black] (10,-2.5) circle [radius=0.1];
\draw [fill=black] (11,0) circle [radius=0.1];
\draw [fill=black] (13,0) circle [radius=0.1];
\draw [fill=black] (15,0) circle [radius=0.1];
\draw [fill=black] (17,0) circle [radius=0.1];

\draw [fill=black] (1,-2) circle [radius=0.1];
\draw [fill=black] (3,-2) circle [radius=0.1];
\draw [fill=black] (5,-2) circle [radius=0.1];
\draw [fill=black] (7,-2) circle [radius=0.1];
\draw [fill=black] (10,-4.5) circle [radius=0.1];
\draw [fill=black] (11,-2) circle [radius=0.1];
\draw [fill=black] (13,-2) circle [radius=0.1];
\draw [fill=black] (15,-2) circle [radius=0.1];
\draw [fill=black] (17,-2) circle [radius=0.1];
}
\end{center}
\caption{A cycle extension.}
\label{fig:cycleextension2}
\end{figure}

{\bf Extension-Rotation algorithm}\\
The algorithm begins by choosing an arbirtrary edge $e\in M$ and letting path $P_1=e$. 
\begin{description}
\item[Basic Idea] 
It proceeds in rounds. At the beginning of round $k$ we will have a 2-path $P_k=(x_1,x_2,\ldots,x_{2k})$. A round consists of the following: Let $Q_0=P_k$ and then for $i=1,2,\ldots$, if necessary, grow a set of paths $Q_1,Q_2,\ldots$. Each $Q_i$ is obtained from some $Q_j,j<i$ by a single rotation. 

We continue until either we make a simple extension (defined below) or we make a cycle extension (defined below) or fail.
\item[Simple Extensions]
The process is curtailed if at any point the procedure generates a path $P=(y_1,y_2,\ldots,y_{2k})$ and an edge $\set{u,v}\in M$ disjoint from $V(P)$ such that\\
 $(y_1,y_2,\ldots,y_{2k},u,v)$ is a 2-path. In which case we can extend our current 2-path to one of length $2k+2$ and end the round. We call this a {\em simple extension}.
\item[Cycle extensions]
If we do not find a simple extension, then we see if there is a path $P=(y_1,y_2,\ldots,y_{2k})\in \cP_L$ such that $\G$ contains the path $(y_{2k-1},y_1,y_{2k},y_2)$ (see Figure~\ref{fig:cycleextension}).

\begin{figure}[h]
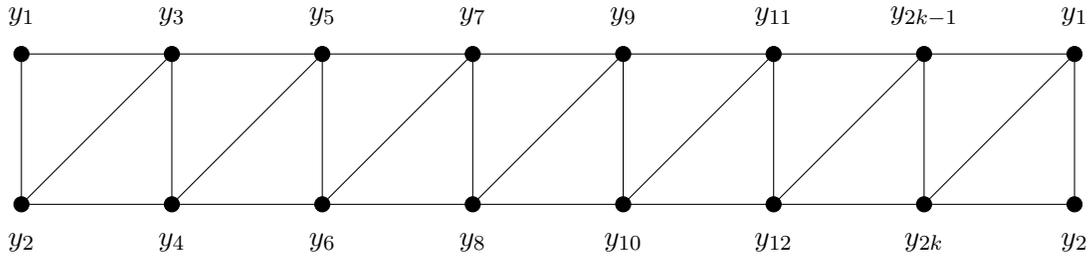

\begin{center}
\pic{
\draw (1,0) -- (15,0);
\draw (1,-2) -- (15,-2);
\draw (1,0) -- (1,-2);
\draw (3,0) -- (3,-2);
\draw (5,0) -- (5,-2);
\draw (7,0) -- (7,-2);
\draw (9,0) -- (9,-2);
\draw (11,0) -- (11,-2);
\draw (13,0) -- (13,-2);
\draw (15,0) -- (15,-2);
\draw (1,-2) -- (3,0);
\draw (3,-2) -- (5,0);
\draw (5,-2) -- (7,0);
\draw (7,-2) -- (9,0);
\draw (9,-2) -- (11,0);
\draw (11,-2) -- (13,0);
\draw (13,-2) -- (15,0);
\node at (1,.5) {$y_1$};
\node at (3,.5) {$y_3$};
\node at (5,.5) {$y_5$};
\node at (7,.5) {$y_7$};
\node at (9,.5) {$y_9$};
\node at (11,.5) {$y_{11}$};
\node at (13,.5) {$y_{2k-1}$};
\node at (1,-2.5) {$y_2$};
\node at (3,-2.5) {$y_4$};
\node at (5,-2.5) {$y_6$};
\node at (7,-2.5) {$y_8$};
\node at (9,-2.5) {$y_{10}$};
\node at (11,-2.5) {$y_{12}$};
\node at (13,-2.5) {$y_{2k}$};
\node at (15,.5) {$y_1$};
\node at (15,-2.5) {$y_2$};

\draw [fill=black] (1,0) circle [radius=0.1];
\draw [fill=black] (3,0) circle [radius=0.1];
\draw [fill=black] (5,0) circle [radius=0.1];
\draw [fill=black] (7,0) circle [radius=0.1];
\draw [fill=black] (9,0) circle [radius=0.1];
\draw [fill=black] (11,0) circle [radius=0.1];
\draw [fill=black] (13,0) circle [radius=0.1];
\draw [fill=black] (15,0) circle [radius=0.1];

\draw [fill=black] (1,-2) circle [radius=0.1];
\draw [fill=black] (3,-2) circle [radius=0.1];
\draw [fill=black] (5,-2) circle [radius=0.1];
\draw [fill=black] (7,-2) circle [radius=0.1];
\draw [fill=black] (9,-2) circle [radius=0.1];
\draw [fill=black] (11,-2) circle [radius=0.1];
\draw [fill=black] (13,-2) circle [radius=0.1];
\draw [fill=black] (15,-2) circle [radius=0.1];
}
\end{center}
\caption{Closing a 2-path with $(y_{2k-1},y_1,y_{2k},y_2)$.}
\label{fig:cycleextension}
\end{figure}

We say that we {\em close the path} to create a cycle $C=(y_1,y_2,\ldots,y_{2k},y_1)$. If $k=n/2$ then we have found the square of a Hamilton cycle. Otherwise, we seek a {\em cycle extension}. By this we mean that find  an edge $\set{u,v}\in M$ disjoint from $V(P)$ such that and $1\leq \ell<k$ such that $G+X$ contains the path  $(y_{2\ell-1},u,y_{2\ell},v)$. In which case we now have the 2-path $(y_{2\ell+1},y_{2\ell+2},\ldots,y_{2k},y_{2k-1},y_1,y_2,\ldots,y_{2\ell-1},y_{2\ell},u,v)$. We call this a {\em cycle extension} (see Figure~\ref{fig:cycleextension2}). If no such pair $\ell,\set{u,v}$ exists then we fail.
\end{description}

We can use all edges of $\G$ at any stage of the algorithm. However, in the description below, we only mention edges that are needed in the analysis. Here we rely on the fact that adding edges will not prevent a successful execution of {\sc ERA}.
\begin{enumerate}[{\bf Step 1}]
\item\label{step:1} We start a round with $P_k=(x_1,x_2,\ldots,x_{2k})$. We then do a set of rotations with $e=\set{x_1,x_2}$ as the fixed end, one for each neighbor of $\set{x_{2k-1},x_{2k}}$ in $\G_1$. Assuming there are no simple extensions we generate a set of 2-paths $Q_1,Q_2,\ldots,Q_L,L\geq \b_1n$. The end pillar of $P_i$, other than $\set{x_1,x_2}$, will be denoted by $e_i$  for $i=1,2,\ldots,L$.  
\item\label{step:2} After this, we take each $Q_i$ in turn and do a set of rotations with $e_i$ as the fixed end and $e$ as the rotated end, using the edges of $\G_2$ for this purpose.  
\item\label{step:3} If we fail to obtain a simple extension, then we use the $\G$ edges to look for a cycle extension, using all of the 2-paths generated for this task.
\end{enumerate}
\subsection{Analysis of {\sc ERA}}
\lem{lem3}{
W.h.p. algorithm {\sc ERA} succeeds in finding the square of a Hamilton cycle.
}{
We argue that w.h.p. we can always find an $X_3$ edge to close a path if there is no simple extension. 
Let $\cP_k$ be the set of 2-paths generated in round $k$ (Step~\ref{step:1} and \ref{step:2}). Thus, 
\[
|\cP_k| \ge L \cdot L/2 \ge \b_1^2n^2/2
\]
and consequently
\[
\Pr(\text{No path of $\cP_k$ can be closed})\leq (1-p^3)^{\b_1^2n^2/2}\leq e^{-\frac12\b_1^2K^3\log n}.
\]
Since there at most $n/2$ rounds we see that w.h.p. there is at least one path in a round that can be closed, if needed.

Having closed a 2-path, the existence of $\ell,\set{u,v}$ follows from the connectivity of $\G_2$, see Lemma~\ref{Gamcon}.
}

\begin{figure}[h]
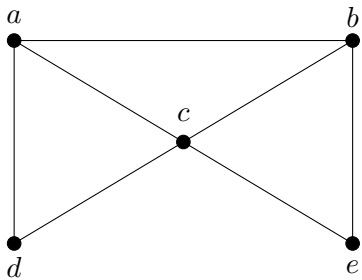

\begin{center}
\pics{0.9}{
\draw (1,0) -- (6,0);
\draw (1,0) -- (1,-3);
\draw (6,0) -- (6,-3);
\draw (1,0) -- (3.5,-1.5);
\draw (6,0) -- (3.5,-1.5);
\draw (1,-3) -- (3.5,-1.5);
\draw (6,-3) -- (3.5,-1.5);
\node at (1,.35) {$a$};
\node at (6,.35) {$b$};
\node at (3.5,-1.1) {$c$};
\node at (1,-3.35) {$d$};
\node at (6,-3.35) {$e$};

\draw [fill=black] (1,0) circle [radius=0.1];
\draw [fill=black] (6,0) circle [radius=0.1];
\draw [fill=black] (1,-3) circle [radius=0.1];
\draw [fill=black] (6,-3) circle [radius=0.1];
\draw [fill=black] (3.5,-1.5) circle [radius=0.1];
}
\end{center}
\caption{A graph used for $n$ odd.}
\label{fig:subgraph}
\end{figure}

When $n=2m+1$ is odd, we use extra $Kn^{4/3}\log^{1/3}n$ edges $X_4$ (chosen independently from $X_1$, $X_2$ and $X_3$) to find a subgraph as in Figure~\ref{fig:subgraph}.

We can then use $a,d$ and $b,e$ as pillars and basically proceed as in the even case, making sure to avoid breaking up this subgraph and follow its vertices as $(d,a,c,b,e)$.

\subsection{A lower bound}
It is as well to consider lower bounds on the number of random edges needs to add to a graph $G$ to obtain the $k$th power of Hamilton cycle. For this we consider the complete bipartite graph $G=K_{s,t}$ with bipartition $A,B$ and where $s=|A|=\a n$ and $s+t=n$. We will be thinking here of the case where $\a$ is a small constant and so it does not fit exactly into the assumptions of Theorem \ref{th1}. In $K_{s,t}$ the lower bound is much less than $n/2$ and have no lower bounds for the case where the minimum degree significantly exceeds $n/2$.

We can associate a sequence $\s_H$ of length $n$ over the alphabet $\set{A,B}$ with a Hamilton cycle $H$ in $K_{s,t}$. The $i$th symbol will be an $A$ if and only if the $i$th vertex of the cycle is in $A$. Only $AB$ edges are in $G$ and it is not difficult to show by examining $\s_H$ that at most $2ks$ of the edges of $H$ can be of this type. It follows that if we add edges to $G$ with probability $p$ then the expected number of $k$th powers will be at most $n!p^{k(n-2s)}$. Thus we require $p\geq n^{-1/k(1-2\a)}$ or at least $n^{2-1/k(1-2\a)}$ random edges. In particular, for $k=2$ this implies that we need $n^{2-1/(2-4\a)}$ random edges,  which for small $\a$ yields $n^{3/2-O(\a)}$. This is close to optimal, since $n^{3/2+o(1)}$ is the trivial upper bound.


\begin{thebibliography}{99}

\bibitem{AKS} M.~Ajtai, J.~Koml\'{o}s and E.~Szemer\'{e}di,
{\em The first occurrence of Hamilton cycles in random graphs},
Annals of Discrete Mathematics~\textbf{27} (1985), 173--178.

\bibitem{BFM} T. Bohman, A.M. Frieze and R. Martin, {\em How many random edges make a dense graph Hamiltonian?}, Random Structures and Algorithms 22 (2003) 33-42.

\bibitem{Boll} B.~Bollob\'{a}s,
{\em The evolution of sparse graphs},
in Graph Theory and Combinatorics, Academic Press, Proceedings of Cambridge Combinatorics,
Conference in Honour of Paul Erd\H{o}s (B.~Bollob\'{a}s; Ed) (1984), 35--57.


\bibitem{ER1} P. Erd\H{o}s and A. R\'enyi, {\em On the existence of a factor of degree one of a connected random graph}, Acta. Math. Acad. Sci. Hungar. 17 (1966) 359-368.


\bibitem{FP} A.M. Frieze and B. Pittel, Perfect matchings in random graphs with prescribed minimal degree, {\em Trends in Mathematics, Birkhauser Verlag, Basel} (2004) 95-132.


\bibitem{KS} J.~Koml\'os and E.~Szemer\'edi,
{\em Limit distributions for the existence of Hamilton circuits in a random graph}, Discrete Mathematics~\textbf{43} (1983), 55--63.

\bibitem{KSS} J.~Koml\'os, N. Sark\"ozy and E.~Szemer\'edi, {\em Proof of the Seymour conjecture for large graphs}, Annals of 

\bibitem{KO} D. K\"uhn and D. Osthus, 
{\em On P\'osa's conjecture for random graphs}, SIAM Journal on Discrete Mathematics \textbf{26} (2012), 1440--1457.

\bibitem{McD} C. McDiarmid, On the method of bounded differences, {\em in Surveys in Combinatorics, ed. J. Siemons, London Mathematical Society Lecture Notes Series} 141, Cambridge University Press, 1989.

\bibitem{NS} R. Nenadov and N. \v{S}kori\'c, 
{\em Powers of cycles in random graphs and hypergraphs}, arxiv:1601.04034v1. 


\bibitem{Posa} L. P\'osa, 
{\em Hamiltonian circuits in random graphs}, Discrete Mathematics~\textbf{14} (1976), 359--364.

\bibitem{Rio} O. Riordan, 
{\em Spanning subgraphs of random graphs}, Combinatorics, Probability and Computing~\textbf{9} (2000), 125--148.

\end{thebibliography}
\end{document}